\documentclass[10pt]{amsart}
\usepackage{amssymb,latexsym, amsthm,epsf}
\newtheorem{thm}{Theorem}[section] 
\newtheorem{defn}[thm]{Definition}

\newtheorem{lem}[thm]{Lemma}

\newtheorem{prop}[thm]{Proposition}
\newtheorem{rem}[thm]{Remark}

\def\N{{\mathbb N}}
\def\R{{\mathbb R}}

\def\({\left(}
\def\){\right)}

\long\def\forget#1\forgotten{}

\newcommand{\sumlim}{\sum\limits}

\newif \iffurther 
\furtherfalse

\newif \iffurther 
\furtherfalse

\newif\ifXY 
\XYfalse    

\ifXY
\usepackage{xy}
\fi
\ifXY
\xyoption{all}
\fi

\begin{document}

\title{The Orchard crossing number of an abstract graph}

\author{Elie Feder and David Garber}

\address{Kingsborough Community College of CUNY, Department of Mathematics and Computer Science,
2001 Oriental Blvd., Brooklyn, NY 11235, USA}
\email{eliefeder@gmail.com, efeder@kbcc.cuny.edu}

\address{Department of Applied Mathematics, Holon Institute of Technology, Golomb 52,
PO Box 305, Holon 58102, Israel}
\email{garber@hit.ac.il}

\date{\today}

\begin{abstract}
We introduce the Orchard crossing number, which is defined in
a similar way to the well-known rectilinear crossing number. We
compute the Orchard crossing number for some simple families of
graphs. We also prove some properties of this crossing number.

Moreover, we define a variant of this crossing number which is
tightly connected to the rectilinear crossing number, and compute it
for some simple families of graphs.
\end{abstract}

\maketitle

\section{Introduction}
Let $G$ be an abstract graph. Motivated by the Orchard relation,
introduced in \cite{bacher,BaGa}, one can define the {\it Orchard
crossing number} of $G$, in a similar way to the well-known {\it
rectilinear crossing number} of an abstract graph $G$ (denoted by
$\overline{\rm cr}(G)$, see \cite{aak2,PT}).

The Orchard crossing number is interesting for several reasons.
First, it is based on the Orchard relation which is an equivalence
relation on the vertices of a graph, with at most two equivalence
classes (see \cite{bacher}). Moreover, since the Orchard relation
can be defined for higher dimensions too (see \cite{bacher}),
hence the Orchard crossing number
may be also generalized to higher dimensions.

Second, a variant of this crossing
number is tightly connected to the well-known rectilinear crossing number (see
Proposition \ref{max_complete}).

Third, one can find real problems which the Orchard crossing number
can represent.  For example, design a network of computers which should be
constructed in a manner which allows possible extensions
of the network in the future. Since we want to avoid (even future) crossings
of the cables which are connecting between the computers, we need to count not
only the present crossings, but also the separators
(which might come to cross in the future).

In the current paper, we define the Orchard crossing number for an
abstract graph $G$, and we compute it for some simple families of graphs.
We also deal with some simple properties of this crossing
number. Some more properties can be found at \cite{GPS}.

We also define some variant of this crossing number, which is
tightly connected to the rectilinear crossing number, and we compute
it for some simple families of graphs.

\medskip

The paper is organized as follows. In Section \ref{defs}, we present
the Orchard relation, define the Orchard crossing number and its
variant, and give some examples.
In Section \ref{crossing_complete}, we compute the Orchard crossing number for
the complete graph $K_n$.
Section \ref{star} deals with the Orchard crossing number of the star graph $K_{n,1}$.
In Section \ref{wheel}, we compute the Orchard crossing number for
the wheel graph $W_{n,1}$.
Section \ref{bipartite} deals with some partial results about complete
bipartite graphs $K_{n,m}$.
In Section \ref{union}, we discuss the computation of the Orchard crossing number 
for a particular case of a union of two graphs.

\section{The Orchard relation and the Orchard crossing
numbers}\label{defs}

We start with some notations.
A finite set $\mathcal P = \{P_1,\cdots, P_n\}$ of $n$ points
in the plane $\R^2$ is a {\it generic configuration}
if no three points of $\mathcal P$ are collinear.

A line $L \subset \R^2$ {\it separates} two points $P,Q \in
(\R^2\setminus L)$ if $P$ and $Q$ are in different connected
components of $\R^2 \setminus L$. Given a generic configuration
$\mathcal P$, denote by $n(P,Q)$ the number of lines defined by
pairs of points in $\mathcal P \setminus \{P,Q\}$, which separate
$P$ and $Q$.

In this situation, one can define:
\begin{defn}[Orchard relation]
For two distinct points $P,Q$ of a generic configuration ${\mathcal P}$,
we set $P \sim Q$ if we have
$$n(P,Q) \equiv (n-3) \pmod 2.$$
\end{defn}

One of the main results of \cite{bacher} is that this relation is an
equivalence relation, having at most two equivalence classes.
Moreover, this relation can be used as an (incomplete) distinguishing
invariant between generic configurations of points.

\medskip

In order to define the {\it Orchard crossing number}, we need some more notions.

\begin{defn}[Rectilinear drawing of an abstract graph $G$]
Let $G=(V,E)$ be an abstract graph.
{\em A rectilinear drawing of a graph $G$}, denoted by $R(G)$, is a generic subset of points $V'$
in the affine plane, in bijection with $V$. An edge $(s,t)\in E$ is
represented by the straight segment $[s',t']$ in $\R^2$.
\end{defn}

One then can associate a {\it crossing number} to such a drawing:
\begin{defn}
Let $R(G)$ be a rectilinear drawing of the abstract graph $G=(V,E)$.
The {\em crossing number} of $R(G)$, denoted by $n(R(G))$, is:
$$n(R(G)) = \sum _{(s,t) \in E} n(s,t)$$
\end{defn}

Note that the sum is taken only over the edges of the graph, whence $n(s,t)$
counted in {\it all} the lines generated by pairs of points of the configuration.

\medskip

Now, we can define the {\it Orchard crossing number} of an abstract
graph $G=(V,E)$:

\begin{defn}[Orchard crossing number]
Let $G=(V,E)$ be an abstract graph.
The {\em Orchard crossing number} of $G$, ${\rm OCN}(G)$, is
$${\rm OCN}(G) = \min _{R(G)} (n(R(G)))$$
\end{defn}

\medskip

Note that the usual intersection of two edges contributes two separations:
Each edge separates the two
points at the ends of the other edge. Hence, one have the following
important observation:
\begin{rem}
Let $G$ be an abstract graph. Then:
$$\overline{\rm cr}(G) \leq \frac{1}{2} {\rm OCN}(G)$$
where the rectilinear crossing number of an abstract graph $\overline{\rm cr}(G)$
is defined to be the smallest number of crossings between edges in a
rectilinear drawing of the graph $G$.
\end{rem}

Hence, by computing the Orchard crossing number of a graph, one gets
also some upper bounds for the rectilinear crossing number (these
bounds are rather bad, see the table after Proposition
\ref{complete}. The reason for the bad bounds is that the best drawing
with respect to the Orchard crossing number is rather worse with
respect to the rectilinear crossing number, see Proposition
\ref{complete}).

\medskip

A variant of the Orchard crossing number is the {\it maximal Orchard crossing number}:

\begin{defn}[Maximal Orchard crossing number]
Let $G=(V,E)$ be an abstract graph. The {\em maximal Orchard crossing
number} of $G$, ${\rm MOCN}(G)$, is
$${\rm MOCN}(G) = \max _{R(G)} (n(R(G)))$$
\end{defn}

This variant is extremely interesting due to the following result:

\begin{prop}\label{max_complete}
The drawing which yields the maximal Orchard crossing number for
complete graphs $K_n$ is the same as the drawing which attains the
rectilinear crossing number of $K_n$.
\end{prop}

The importance of this result is that it might be possible that the
computation of the maximal Orchard crossing number is easier than
the computation of the rectilinear crossing number.

\begin{proof}
Let us look on quadruples of points: there are only two
possibilities to draw four points in a generic position, see Figure \ref{quadraple}.

\begin{figure}[h]
\epsfysize=3cm \centerline{\epsfbox{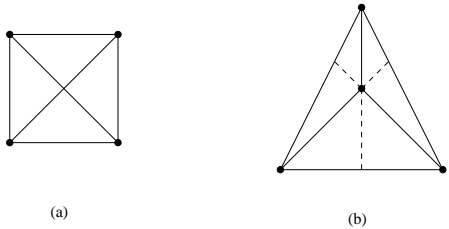}} \caption{Two
generic configurations of $4$ points} \label{quadraple}
\end{figure}

Note that $\overline{\rm cr}(K_4)=0$, since we can draw $K_4$ without
any crossings of the edges (type (b)). Hence, for getting the
minimal rectilinear drawing for $K_n$ with respect to the
rectilinear crossing number, we would like to have as many as
possible quadruples of points arranged as in type (b).

On the other hand, for getting the maximal number of crossings with
respect to the Orchard crossing number, we would like again to have
as many as possible quadruples of points arranged as in type (b),
since type (a) contributes 2 crossings (each internal edge
separates the two ends of the other internal edge) while type (b)
contributes 3 crossings (one on each edge of the convex hull).
Hence, the rectilinear drawing which yields the maximal number of
crossings with respect to the Orchard crossing number, will also give
the minimal number of crossings with respect to the rectilinear
crossing number.
\end{proof}

\subsection{Some examples}

In this subsection, we give some examples for computing the Orchard
crossing numbers for some small graphs.

\medskip

The Orchard crossing number of $K_4$ is $2$, since there are only
two generic drawings of four points as presented in Figure
\ref{quadraple}. Type (a) has $2$ crossings, while type (b) has $3$
crossings. Therefore:
$${\rm OCN}(K_4)=2.$$
$${\rm MOCN}(K_4)=3.$$

\medskip

\begin{figure}[h]
\epsfysize=2.5cm \centerline{\epsfbox{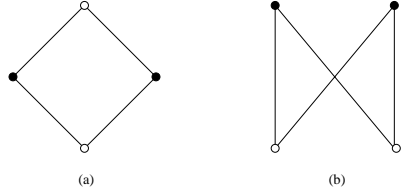}} \caption{Generic rectilinear drawings of $K_{2,2}$}\label{4bipartite}
\end{figure}

For the complete bipartite graph $K_{2,2}$, we have more
possibilities for rectilinear drawings, since except for the two
choices for drawings, we have also to choose the two pairs of
points. Figure \ref{4bipartite}(a) shows a drawing of $K_{2,2}$ without
crossings at all (we distinguish between the two sets of two points
by their colors). On the other hand,  Figure \ref{4bipartite}(b) shows a
rectilinear drawing of $K_{2,2}$ with $2$ crossings. Hence, we get that
$${\rm OCN}(K_{2,2})=0.$$
$${\rm MOCN}(K_{2,2})=2.$$

Note that $K_{2,2} \cong C_4$ where $C_4$ is the cycle graph on $4$
vertices, and it is shown in \cite{GPS} that the graphs whose
Orchard crossing number equals $0$ are the following:
\begin{enumerate}
\item The cycle graphs, $C_n$, and their subgraphs.
\item The graph on $4$ vertices presented in Figure \ref{graph_4}.

\begin{figure}[h]
\epsfysize=2cm \centerline{\epsfbox{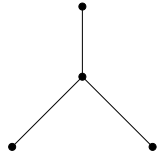}} \caption{An example of a graph
with ${\rm OCN}(G)=0$} \label{graph_4}
\end{figure}

\end{enumerate}

\section{Orchard crossing numbers for complete graphs $K_n$}\label{crossing_complete}

For the complete graphs $K_n$, we have the following proposition:

\begin{prop}\label{complete}
Let $n \in \N$. Then
$${\rm OCN}(K_n)=2{n \choose 4},$$
which is obtained by placing all vertices of $K_n$ in a convex
position.
\end{prop}

\begin{proof}
Since there are only two generic rectilinear drawings of four points
(as presented in Figure \ref{quadraple}), any quadruple of points
can contribute $2$ (type (a)) or $3$ (type (b)) to the total number
of crossings.

Since we want to minimize the total number of crossings, we have to
require that every quadruple of  points will be in a convex
position. This immediately implies that all $n$ points have to be in
 a convex position  as needed.

The computation of the crossing number is straight-forward: from the
configuration of $n$ points, we can choose ${n \choose 4}$ different
quadruples of points, and each quadruple contributes $2$ to the
total number of crossings, which implies that the total number of
crossings is $2{n \choose 4}$.
\end{proof}

Based on the last proposition, we compare in the following table between
the Orchard crossing numbers and the rectilinear crossing numbers
for the complete graph $K_n$, for $n \leq 12$.

\medskip

$$\begin{array}{|c|c|c|}
\hline
n   & {\rm OCN}(K_n)         & \overline{\rm cr}(K_n)      \\
\hline
 4  & 2                      &  0                     \\
 5  & 10                     &  1                     \\
 6  & 30                     &  3                     \\
 7  & 70                     &  9                     \\
 8  & 140                    &  19                    \\
 9  & 252                    &  36                    \\
 10 & 420                    &  62                    \\
 11 & 660                    &  102                   \\
 12 & 990                    &  153                   \\
\hline
\end{array}$$

\medskip

Note that up to $n = 27$, $\overline{\rm cr}(K_n)$ is known, and for $n \geq 28$,
only bounds are known (see \cite{AFLS,aak2,AGOR,ak2,aor}).

\section{Orchard crossing numbers for the star graph $K_{n,1}$}\label{star}

In this section, we deal with the Orchard crossing number of the graph $K_{n,1}$ (which is also called the {\it star graph}).
Assume that we have $n$ black points and one white point, and the white point is connected to all the black points.
Then, the following holds:

\begin{prop}\label{k_m1_prop}
The configuration of $n+1$ points which attains the minimal number of
Orchard crossings  for the graph $K_{n,1}$ consists of $n$ black points which are the
vertices of a regular $n$-gon, and the white point is located at the
center of this polygon for odd $n$, or close to its center (i.e. it
is not separated by a line from the center), for even $n$ (see
Figure \ref{k_m1} for examples of $n=7,8$).

\begin{figure}[h]
\epsfysize=5cm \centerline{\epsfbox{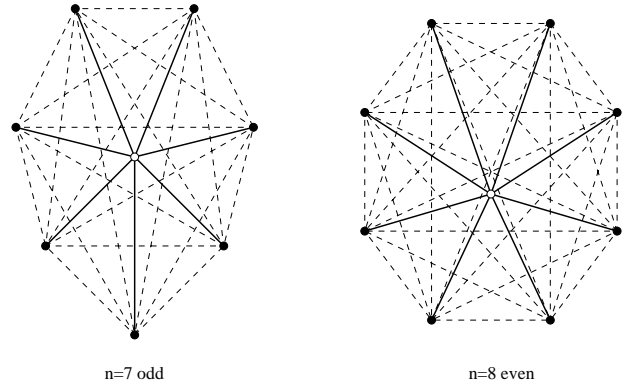}} \caption{Examples of
$K_{n,1}$ for odd $n$ and for even $n$} \label{k_m1}
\end{figure}

Hence, we have:
$$
{\rm OCN}(K_{n,1})=\left\{
\begin{array}{cl}
            \frac{n(n-1)(n-3)}{8}  &  {\rm odd}\ n \\
            \frac{n(n-2)^2}{8}     &  {\rm even}\ n
\end{array}
\right.
$$

\end{prop}

\begin{proof}
As in the proof of Proposition \ref{complete},
we start by looking on quadruples of points. We have three cases of
three black points and one white point (we exclude the case of quadruples of
four black points, since it contributes nothing in any position):
\begin{enumerate}
\item A white point inside a triangle of black points
(see Figure \ref{bipartite3-1}(a)). This quadruple contributes
nothing to the total number of crossings.
\item Three black points and one white point in a convex position
(see Figure \ref{bipartite3-1}(b)). This quadruple contributes $1$
to the total number of crossings.
\item A black point inside a triangle consists of one white point and
two black points (see Figure \ref{bipartite3-1}(c)). This quadruple
contributes $2$ to the total number of crossings.
\end{enumerate}

\begin{figure}[h]
\epsfysize=3cm \centerline{\epsfbox{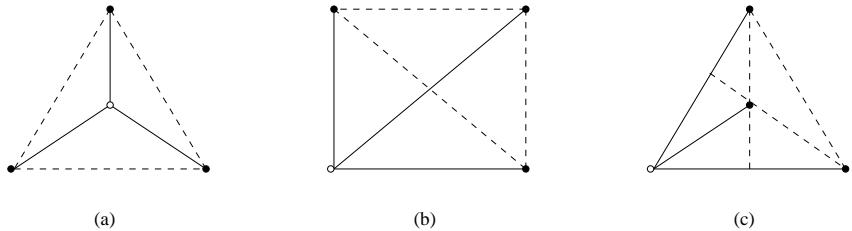}}
\caption{Three types of quadruples of points with three black points
and one white point} \label{bipartite3-1}
\end{figure}

Our aim is to minimize the number of quadruples of types (2) and (3)
as a whole. Later, we will minimize the number of quadruples of type
(3) as opposed to those of type (2). In order to compute the number
of quadruples of types (2) and (3) in a given drawing of $K_{n,1}$,
we label the black points $b_1, b_2,\dots,b_n$ and the white point
$w$.

For each $1 \leq i \leq n$, let $L_i$ be the line which connects
$b_i$ to $w$. This line divides the remaining $n-1$ black points
into two subsets - one contains $k_i$ points, and the other
contains $(n-1)-k_i$ points. For each pair of points, $b_j$ and
$b_k$, situated on the same side of $L_i$, the quadruple $b_i, b_j,
b_k, w$, is of type (2) or (3) and hence contributes at least one
crossing. We have ${k_i \choose 2} + {(n-1)-k_i \choose 2}$ such
pairs of points. Thus in total we have $\sumlim_{i=1}^n \left( {k_i
\choose 2} + {(n-1)-k_i \choose 2} \right)$ quadruples of types (2)
or (3). (Note that if $b_j$ and $b_k$ are separated by $L_i$, then
the quadruple $b_i, b_j, b_k, w$, is of type (1) or (2) and hence
contributes at most one crossing.)

Note that each quadruple is counted twice due to the following
reason. Consider the quadruple $b_i, b_j, b_k, w$. Assume, without
loss of generality, that $b_i, b_j$ and $b_k$ are oriented in
clockwise fashion around $w$. Then, both lines $L_i$ and $L_k$ give
this quadruple. However, the line $L_j$ separates $b_i$ from $b_k$
and therefore will not give this quadruple. So we divide the above
number of quadruples by $2$.

Since each quadruple of type (2) or type (3) contributes at least one
crossing, we have at least $$\frac{\sumlim_{i=1}^n \left( {k_i
\choose 2} + {(n-1)-k_i \choose 2} \right) }{2}$$ crossings.

This sum is minimized if $k_i = \frac{n-1}{2}$ for all $i$, i.e.
every line connecting the white point to a black point divides the
remaining black points into two subsets with the same number of
points (or different by 1 for even $n$). This condition will be
satisfied only if the $n$ points are evenly distributed around the
white point.

In order to determine the distance of the $n$ black points to the white
point, we consider types (2) and (3). In order to avoid any
quadruple of type (3), we have to ensure that each quadruple $b_i,
b_j, b_k, w$ is in a convex position. This will be accomplished only
if the $n$ black points are all in a convex position around the white point.
In this case, every quadruple which contains a white point, will be
of type (2) and therefore will contribute only one crossing. This
gives us the requested minimal rectilinear drawing.

\medskip

Therefore, we can compute the number of crossings as follows. For
odd $n$, we have $k_i=\frac{n-1}{2}$ for all $1 \leq i \leq n$.
Thus, we have:
\begin{eqnarray*}
{\rm OCN}(K_{n,1}) &= & \frac{\sumlim_{i=1}^n \left( {\frac{n-1}{2}
\choose 2} +
{(n-1) - \frac{n-1}{2} \choose 2} \right) }{2} = \\
&= & n {\frac{n-1}{2} \choose 2} = \frac{n\left(\frac{n-1}{2}\right)\left(\frac{n-3}{2}\right)}{2} =\\
& =& \frac{n(n-1)(n-3)}{8}.
\end{eqnarray*}

For even $n$, we have $k_i=\frac{n}{2}$ for all $1 \leq i \leq n$.
Thus, we have:
\begin{eqnarray*}
{\rm OCN}(K_{n,1}) & =& \frac{\sumlim_{i=1}^n \left( {\frac{n}{2}
\choose 2} +
{(n-1) - \frac{n}{2} \choose 2}\right)}{2} =\\
&= & \frac{\sumlim_{i=1}^n \left( {\frac{n}{2} \choose 2} + { \frac{n-2}{2} \choose 2} \right)}{2} = \\
&=& n \cdot \frac{\frac{\left(\frac{n}{2}\right)
\left(\frac{n-2}{2}\right)}{2} +
\frac{\left(\frac{n-2}{2}\right) \left(\frac{n-4}{2}\right) }{2}}{2}=\\
&=& \frac{n(n-2)(n+n-4)}{16}=\\
&=& \frac{n(n-2)(2n-4)}{16} = \frac{n(n-2)^2}{8}.
\end{eqnarray*}

\end{proof}

In the following lemma, we present a different way to count the number of Orchard
crossings in the minimal rectilinear drawing, as in Proposition
\ref{k_m1_prop}.

\begin{lem}
$$
{\rm OCN}(K_{n,1})=\left\{
\begin{array}{cl}
            \frac{n(n-1)(n-3)}{8}  &  {\rm odd}\ n \\
            \frac{n(n-2)^2}{8}     &  {\rm even}\ n
\end{array}
\right.
$$
\end{lem}

\begin{proof}
Since there are no quadruples of type (3) (since the only internal
point is white), all the crossings come from quadruples of type (2),
which contribute one crossing for each quadruple. So, we have to
count the number of quadruples of type (2): quadruples with three
black points and one white point in a convex position.

In the odd case, since the white point is in the center of the
$n$-gon, for having a quadruple of type (2), all the black points should
be within $180^{\circ}$ of each other. For counting the number of
such triples, we fix one black point. The other two points should be
within $180^{\circ}$ from this point. This gives us $\frac{n-1}{2}$
options to choose two points from. We can do this in ${\frac{n-1}{2}
\choose 2}$ ways. Now since we could have started with any of the
$n$ points, we have to multiply it by $n$. So, we get:
$$n {\frac{n-1}{2} \choose 2} = \frac{n \cdot \frac{n-1}{2}
\frac{n-3}{2}}{2} = \frac{n(n-1)(n-3)}{8},$$ as needed.

\medskip

Now we turn to the even case. This case is more complicated, since
if we choose two antipodal black points, we have to cut the number
of quadruples by $2$, due to the central white point which is not
located exactly at the center (see Figure \ref{k_m1}, $n=8$). So, we
split the count of triples into two cases: triples which do not
include antipodal points, and triples which include antipodal
points.

We start with the case which does not include antipodal points: Fix a
black point. We have to choose two black points out of $\frac{n-2}{2}$
options. Hence,  we have ${\frac{n-2}{2} \choose 2}$ possibilities , and we have in total:
$$n {\frac{n-2}{2} \choose 2} = \frac{n(n-2)(n-4)}{8}$$   crossings.

Now we count the number of triples which include antipodal points:
We have $\frac{n}{2}$ pairs of points in antipodal position.
The third point can be any of the other $n-2$ points. Since
we only count half of these triples (as explained above), we have
$\frac{n(n-2)}{4}$ crossings.

Summing up the two cases, we get:
$$ \frac{n(n-2)(n-4)}{8} + \frac{n(n-2)}{4} = \frac{n(n-2)(n-2)}{8},$$
again as needed.
\end{proof}

\section{Orchard crossing number of a wheel}\label{wheel}

\begin{figure}[h]
\epsfysize=3cm
\centerline{\epsfbox{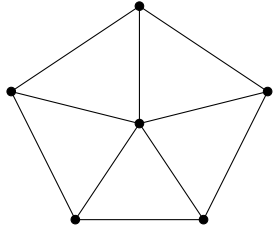}}
\caption{The wheel graph $W_{5,1}$} \label{wheel_fig}
\end{figure}

In this section, we deal with the Orchard crossing number of a {\it
wheel}, which is an abstract graph on $n+1$ points, which is
composed of a cycle of $n$ points and one other point connected
to each point of the cycle (see Figure \ref{wheel_fig} for an example). We denote
such a graph by $W_{n,1}$

Before stating the general case, a simple observation yields the following facts:
\begin{rem}
$${\rm OCN}(W_{3,1})=2$$
$${\rm OCN}(W_{4,1})=6$$
\end{rem}

\begin{proof}
Note that $W_{3,1}=K_4$, and hence the result follows.

For $W_{4,1}$, we have two optimal drawings, each of them has $6$ crossings, see Figure \ref{OCN_draws_W4-1}.

\begin{figure}[h]
\epsfysize=3cm
\centerline{\epsfbox{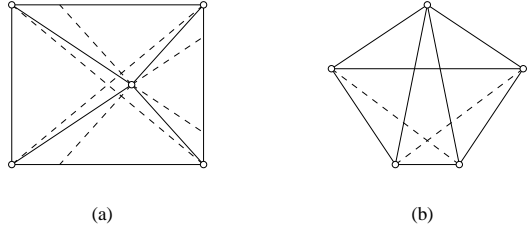}}
\caption{Two optimal drawings for $W_{4,1}$} \label{OCN_draws_W4-1}
\end{figure}

\end{proof}

Here we present the general case:
\begin{prop}\label{wheel_prop}
Let $n \geq 5$. Let $G=W_{n,1}$ be a graph consisting of $n+1$
vertices configured as a wheel.

The configuration of $n+1$ points which attains the minimal number of
Orchard crossings consists of $n$ points which are the vertices of a
regular $n$-gon (these points form the cycle), and the other point
is located at the center of this polygon for odd $n$, or close to
its center (i.e. it is not separated by a line from the center) for
even $n$.

Hence, we have:
$$
{\rm OCN}(W_{n,1})=\left\{
\begin{array}{cl}
            \frac{n(n^2 - 4n + 11)}{8}  &  {\rm odd}\ n \\
            \frac{n(n^2 - 4n + 12)}{8}  &  {\rm even}\ n
\end{array}
\right.
$$

\end{prop}

First we will prove that this drawing minimizes the crossings for a
wheel. Notice that a wheel $W_{n,1}$ is very similar to $K_{n,1}$.
The only difference is that in $W_{n,1}$, we have an additional cycle on the $n$
external points.

For $K_{n,1}$, we know that the above configuration minimizes the
number of crossings. Thus, all we have to show is that the addition
of the cycle maintains the minimality. Hence, we have to consider
how many crossings are added by the addition of the cycle in the
optimal configuration, and how many are added in other
configurations.

There are two types of edges in our configuration:
\begin{enumerate}
\item the edges on the cycle, and

\item the edges which connect the center and an external point.
\end{enumerate}

Let us consider them one at a time.

\begin{enumerate}
\item The edges of the cycle do not separate between any pair of points in the
configuration.

\item Each edge connecting the center and an external point will
intersect the cycle once (where the continuation of the edge crosses
the convex hull). This intersection contributes $1$ to the total
number of crossings (because this edge separates the vertices on the
cycle).
\end{enumerate}

Since we have $n$ external points, we conclude that the addition
of the cycle contributes $n$ crossings in total. Now let us consider
how many crossings the addition of the cycle contributes in other
configurations.

\begin{lem}
Let $R(W_{n,1})$ be a rectilinear drawing for the wheel
$W_{n,1}$. Then, we have for $n \geq 5$:
$$n(R(W_{n,1})) \geq {\rm OCN}(K_{n,1}) + n,$$
\end{lem}

\begin{proof}
We split the proof into two cases, depending on the position
of the central point of the wheel in the drawing.

\medskip

\noindent {\bf Case 1:} If the central point of the wheel is inside
the convex hull generated by the points of the cycle,
then for any edge connecting the center and an
external point, the extension of this edge separates two points on
the convex hull of the cycle. This occurs where the extension of the edge exits the
interior of the cycle. Since we have $n$ points on the cycle, we add at
least $n$ crossings. Hence, we have:
$$n(R(W_{n,1})) \geq n(R(K_{n,1}))+n \geq {\rm OCN}(K_{n,1}) + n,$$
(where $R(K_{n,1})$ is the corresponding drawing for $K_{n,1}$).

\medskip

\noindent {\bf Case 2:} If the central point of the wheel is outside
the convex hull generated by the points of the cycle, we just show that
$$n(R(W_{n,1})) \geq {\rm OCN}(K_{n,1}) + n.$$
by a direct computation: in case that the central point of the wheel is outside
the convex hull generated by the points of the cycle, each quadruple involved the central point and
three points of the cycles contributes at least one crossing (see types (b) and (c) in Figure \ref{bipartite3-1}). Since we have $n \choose 3$ such quadruples,
we have $n(R(W_{n,1})) \geq {n \choose 3}$. On the other hand, for $n \geq 5$,
$${n \choose 3} \geq {\rm OCN}(K_{n,1}) + n$$
since for odd $n$ (the even case is similar), we have:
\begin{eqnarray*}
& &{n \choose 3} - ({\rm OCN}(K_{n,1})\ + n)=\\
&=&\frac{n(n-1)(n-2)}{6}-\frac{n(n-1)(n-3)}{8}-n=\\
&=&\frac{4n(n-1)(n-2)-3n(n-1)(n-3)-24n}{24}=\\
&=&\frac{n^3-25n}{24}\geq 0
\end{eqnarray*}
\end{proof}

Hence, we conclude that the rectilinear drawing of the wheel with the minimal
number of Orchard crossings is as stated in Proposition
\ref{wheel_prop}.

To verify that the values of ${\rm OCN}(W_{n,1})$
which stated in Proposition \ref{wheel_prop} are correct, we simply
have to take the values for ${\rm OCN}(K_{n,1})$ and add $n$ for the
$n$ crossings contributed by the edges of the cycle (as above).

\section{Complete bipartite graphs of type $K_{n,m}$}\label{bipartite}

The case of complete bipartite graphs of type $K_{n,n}$ is complicated,
and we leave its long proof of the following proposition to a different paper \cite{FG}:
\begin{prop}
$${\rm OCN}(K_{n,n})=4n{n \choose 3},$$
which is achieved by the regular $2n$-gon, alternating in colors.
\end{prop}

\medskip

In Table \ref{table1}, we summarize our computational results (based on the database of Aichholzer, Aurenhammer and Krasser \cite{aak})
for the Orchard crossing number for some complete bipartite graphs
$K_{n,m}$, where $n \not = m$, $n >1$ and  $m >1$. In this table,
the rows will be the different values of $n$, and the columns will
be the different values of $m$.

\begin{table}[!h]
\centerline{\begin{tabular}{|l|r|r|r|r|r|}
\hline
n:  /  m:           & 2  &  3  & 4  & 5  & 6  \\
\hline
   2                & 0  &  4  & 12 & 26 & 48 \\
   3                & 4  &  12 & 32 & 63 &    \\
   4                & 12 &  32 & 64 &    &    \\
   5                & 26 &  63 &    &    &    \\
   6                & 48 &     &    &    &    \\
\hline
\end{tabular}}
\caption{${\rm OCN}(K_{n,m})$ for different values of $n$ and $m$}\label{table1}
\end{table}

\subsection{The maximal Orchard crossing number for complete bipartite graphs
$K_{n,m}$}\label{max_crossing_bipartite}

In this section, we deal with the maximal Orchard crossing number of bipartite graphs $K_{n,m}$.
Assume that we have $n$ black points and $m$ white points, and all the white points are connected to all the black points.
Then, the following holds:

\begin{prop}\label{mocn_bipartite}
The configuration of $n+m$ points which attains the maximal number of
Orchard crossings  for the graph $K_{n,m}$ consists of $n$ black points and $m$ white points which are organized
on two arcs facing each other, where the black points are located on one arc
and the white points are located on the other arc (see
Figure \ref{mocn_K44} for an example for $n=m=4$).

\begin{figure}[h]
\epsfysize=3cm
\centerline{\epsfbox{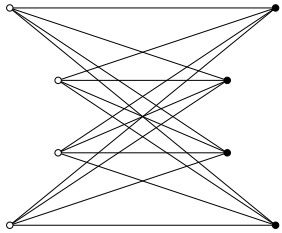}}
\caption{The optimal drawing for ${\rm MOCN}(K_{4,4})$} \label{mocn_K44}
\end{figure}

Hence, we have:
$$
{\rm MOCN}(K_{n,m})=2{n \choose 2}{m \choose 2} + 2m{n \choose 3}+2n{m \choose 3}
$$

\end{prop}

\begin{proof}
As in the previous proofs, we start by looking on quadruples of points. We have three types of quadruples with
two black points and two white points:
\begin{enumerate}
\item Two black points and two white points in a convex position and the points alternate in colors
(see Figure \ref{max_bipartite2-2}(a)). This quadruple contributes
nothing to the total number of crossings.

\item Two black points and two white points in a convex position and the points do not alternate in colors
(see Figure \ref{max_bipartite2-2}(b)). This quadruple contributes
$2$ to the total number of crossings.

\item Two black points and two white points which are not in a convex position
(see Figure \ref{max_bipartite2-2}(c)). This quadruple contributes
$2$ to the total number of crossings.
\end{enumerate}

\begin{figure}[h]
\epsfysize=2.5cm
\centerline{\epsfbox{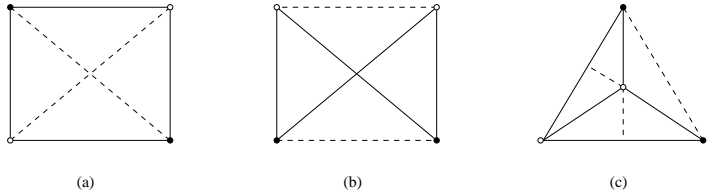}}
\caption{Three cases of quadruples of points with two black points
and two white points} \label{max_bipartite2-2}
\end{figure}

On the other hand, as in the proof of Proposition \ref{k_m1_prop}, we have three types of quadruples with three black points and one white point:
\begin{enumerate}
\item The three black points are in the convex hull and the white point is inside this triangle (see Figure \ref{bipartite3-1}(a)).
This quadruple contributes nothing to the total number of crossings.

\item All the four points are in a convex position
(see Figure \ref{bipartite3-1}(b)). This quadruple contributes
$1$ to the total number of crossings.

\item The convex hull consists of two black points and a white point
(see Figure \ref{bipartite3-1}(c)). This quadruple contributes
$2$ to the total number of crossings.
\end{enumerate}

We shall show that the drawing appearing in the formulation of the proposition indeed attains the maximal possible number
of Orchard crossings for the graph $K_{n,m}$. Let us start with quadruples with two black points and two white points. The maximal number of
crossings which such a quadruple can contribute is $2$. In the above drawing, all the quadruples consisting on two black points
and two white points indeed contribute $2$ (all the quadruples are organized as in type (b)).

Now we move to quadruples with three points of one color and one point of the other color. For such quadruples, the maximal number of
crossings is again $2$. Note that in the above drawing, all the quadruples of this type indeed contribute $2$ (all the quadruples are
organized as in type (c)).

Hence, we get that the drawing appearing in the formulation of the proposition indeed attains the maximal possible number
of Orchard crossings for $K_{n,m}$.

\medskip

For computing ${\rm MOCN}(K_{n,m})$, we simply count the number of quadruples of the different types:
\begin{enumerate}
\item We have ${n \choose 2}{m \choose 2}$ quadruples with two black points and two white points,
since for constructing such a quadruple, we have to choose two black points
(out of $n$) and two white points (out of $m$). Since each quadruple contributes $2$,
the total contribution of these quadruples is $2{n \choose 2}{m \choose 2}$.

\item We have $n{m \choose 3}$ quadruples with one black point and three white points,
since for constructing such a quadruple, we have to choose one black point
(out of $n$) and three white points (out of $m$). Similarly, we have $m{n \choose 3}$
quadruples with one white point and three black points. Since each quadruple contributes $2$,
the total contribution of these quadruples is $2m{n \choose 3}+2n{m \choose 3}$.
\end{enumerate}

Summing up the contributions, we get the desired result.
\end{proof}

\section{The Orchard crossing number of a union of two
graphs}\label{union}

In this section, we deal with the question: what happens to the
Orchard crossing number when we unite two graphs.

We deal with a property of a union of two special graphs:
\begin{prop}
 Let $G$ and $H$ be two graphs on the {\em same} set of vertices.
Assume that $G$ and $H$ get the minimal number of Orchard crossings
in the same rectilinear drawing $C$. Assume also that $H$ is
contained in the complement of $G$ (so they have no edge in common).

Then, the configuration $C$ is the minimal configuration with
respect to the Orchard crossing number for $G \cup H$, and:
$${\rm OCN}(G \cup H)={\rm OCN}(G) + {\rm OCN}(H).$$
\end{prop}

\begin{proof}
Assume on the contrary that there is a better rectilinear drawing
for $G \cup H$. Since in the better rectilinear drawing, ${\rm
OCN}(G \cup H) < {\rm OCN}(G) +{\rm OCN}(H)$, at least one of the
following will occur:

\begin{enumerate}
\item If we delete the edges of $H$, we get a better rectilinear drawing for
$G$, or
\item  If we delete the edges of $G$, we get a better rectilinear
drawing for $H$,
\end{enumerate}
which is a contradiction.
\end{proof}

\begin{rem}
Although one can decompose $W_{n,1}$ into a disjoint union of $K_{n,1}$
and a cycle of length $n$, $C_n$, we can not apply the last proposition, since
the cycle of length $n$ has another isolated vertex, and hence its minimal drawing
is no more $n$ points in convex position and one point inside (as in Figure \ref{Cn_U_pt}(a),
which has $n$ crossings), but $n+1$ points in a convex position (as in Figure \ref{Cn_U_pt}(b), which
has $n-2$ crossings).

\begin{figure}[h]
\epsfysize=3cm
\centerline{\epsfbox{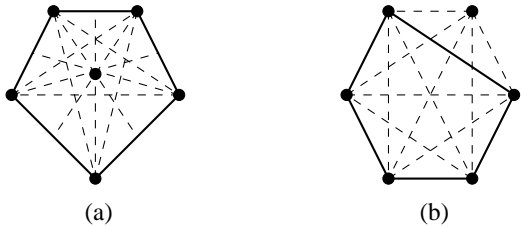}}
\caption{Two drawings of a cycle of length $n=5$ with additional isolated vertex} \label{Cn_U_pt}
\end{figure}

\end{rem}

\section*{Acknowledgments}

We wish to thank Roland Bacher for fruitful discussions. The second
author wish to thank Mikhail Zaidenberg and Institut Fourier for
hosting his stay in Grenoble, where this research was initiated.

\end{document}